\documentstyle [11pt,twoside]   {article}

  \newtheorem{theorem}{Theorem}
  
  \newtheorem{corollary}{Corollary}

  \newcommand {\dz}  {\mbox{\it Proof. \,\,}}

  \newcommand {\qed} {\null \hfill \rule{2mm}{2mm}}

\def  \Z {{\bf Z}}

\def  \Q {{\bf Q}}

\begin {document}

\title{{\large{\bf DIOPHANTINE TRIPLES AND CONSTRUCTION OF HIGH-RANK
       ELLIPTIC CURVES OVER $\Q$ WITH THREE NON-TRIVIAL 2-TORSION POINTS }}}

\author
{\normalsize ANDREJ DUJELLA }

\date{}
\maketitle

\section{Introduction}

Let $E$ be an elliptic curve over $\Q$. The famous theorem of
Mordell-Weil states that
\[ E(\Q) \simeq E(\Q)_{\rm tors} \times {\Z}^r, \]
and by a theorem of Mazur \cite{Maz} we know that only possible torsion
groups over $\Q$ are
\[ E(\Q)_{\rm tors}= \left\{ \begin{array}{ll}
   \Z/m\Z; & \mbox{$m=1,2,\ldots,10$ or $12$}, \\
   \Z/2\Z \times \Z/2m\Z; & \mbox{$m=1,2,3,4$}.
   \end{array}
   \right. \]
Let
\[ B(F)=\sup \{ {\rm rank}\, (E) : \mbox{$E$ curve over $\Q$ with
   $E(\Q)_{\rm tors} \simeq F$} \}, \]
\[ B_r(F)=\lim \sup \{ {\rm rank}\, (E) : \mbox{$E$ curve over $\Q$ with
   $E(\Q)_{\rm tors} \simeq F$} \}. \]
An open question is whether $B(F)<\infty$.

The examples of Martin-McMillen and Fermigier \cite{Fer} show that
$B(0)\geq 23$ and $B(\Z/2\Z)\geq 14$. It follows from results
of Montgomery \cite{Mon} and Atkin-Morain \cite{A-M}
that $B_r(F)\geq 1$ for all torsion groups $F$. Kihara \cite{Kih1}
proved that $B_r(0)\geq 14$ and Fermigier \cite{Fer} that
$B_r(\Z/2\Z)\geq 8$. Recently, Kihara \cite{Kih2} and Kulesz \cite{Kul}
proved using parametrization by $\Q(t)$ and $\Q(t_1,t_2,t_3,t_4)$ that
$B_r(\Z/2\Z \times \Z/2\Z)\geq 4$, and Kihara \cite{Kih3} proved using
parametrization by rational points of an elliptic curve that
$B_r(\Z/2\Z \times \Z/2\Z)\geq 5$. Kulesz also proved that
$B_r(\Z/3\Z)\geq 6$, $B_r(\Z/4\Z)\geq 3$, $B_r(\Z/5\Z)\geq 2$,
$B_r(\Z/6\Z)\geq 2$ and $B_r(\Z/2\Z \times \Z/4\Z)\geq 2$.
The methods used in \cite{Kih2} and \cite{Kul} are similar to the method
of Mestre \cite{Mes1,Mes2}.

In the present paper we prove that $B_r(\Z/2\Z \times \Z/2\Z)\geq 4$
by different method. Namely, we use the theory of, so called, Diophantine
$m$-tuples. By specialization, we obtain an example of elliptic curve
over $\Q$ with torsion group $\Z/2\Z \times \Z/2\Z$ whose rank is equal $7$,
which shows that $B(\Z/2\Z \times \Z/2\Z)\geq 7$.

\section{Construction}
A set of $m$ non-zero rationals $\{a_1,a_2,\ldots,a_m\}$ is called
{\it a (rational) Diophantine $m$-tuple} if $a_ia_j+1$ is a perfect
square for all $1\leq i<j\leq m$ (see \cite{D-acta2}).

Let $\{a,b,c\}$ be a Diophantine triple, i.e.
\[ ab+1=q^2, \quad ac+1=r^2, \quad bc+1=s^2. \]
Define
\[ d=a+b+c+2abc+2qrs, \qquad e=a+b+c+2abc-2qrs. \]
Then it can be easily checked that $ad+1$, $bd+1$, $cd+1$, $ae+1$,
$be+1$ and $ce+1$ are perfect squares. For example,
$ad+1=(as+qr)^2$.

Let us mention that for $a,b,c$ positive integers there is a conjecture that
if $x$ is a positive integer such that $\{a,b,c,x\}$ is a Diophantine
quadruple, then $x$ has to be equal to $d$ or $e$. This conjecture
was verified for some special Diophantine triples (see
\cite{B-D,D-pdeb1,D-pams,D-P,Ked}).

Furthermore, assume that $de+1$ is also a perfect square. Note that this is
impossible if $a,b,c$ are positive integers and $de\neq 0$, but it is
possible for rationals $a,b,c$.

Consider now the elliptic curve
\[ E: \quad y^2=(bx+1)(dx+1)(ex+1). \]
One may expect that $E$ has at least four independent points of
infinite order, namely, points with $x$-coordinates
\[ 0, \quad a, \quad c, \quad \frac{1}{bde} \,. \]

The main problem is to satisfy condition $de+1=w^2$. It can be done,
for example, in the following way. Let $a$ be fixed. Put
$b=ak^2+2k$. Then $q=ak+1$, and put $c=4q(q-a)(b-q)$. It is easy
to check that now $\{a,b,c\}$ is a Diophantine triple. Namely,
$r=q^2+ab-2aq$ and $s=q^2+ab-2bq$. Furthermore, let $ak=t$.
Now the condition $de+1=w^2$ becomes
\begin{equation}  \label{1}
  \begin{array}{c}
  [k^2(t\!+\!2)(2t\!+\!1)(2t\!+\!3)-4k(t\!+\!1)(2t^2\!+\!4t\!+\!1)
   +t(2t\!+\!1)(2t\!+\!3)]^2 \vspace{1ex}\\
  - k^2(4t^2+8t+3)=w^2.
    \end{array}
\end{equation}
There are two obvious solutions of (\ref{1}), namely,
$(k_0,w_0)=(0,t(2t+1)(2t+3))$ and $(k_1,w_1)=(1,1)$, but in both cases we
have $bcd=0$, and therefore they do not lead to a usable formula.
However, using the  solution $(k_0,w_0)$ we may construct
a non-trivial and usable solution of (\ref{1}).
Denote the polynomial on the left side of
(\ref{1}) by $F(k,t)$. Choose the polynomial
$f(k,t)=\alpha(t)k^2+\beta(t)k+\gamma(t)$ such that
\[ F(k,t)-[f(k,t)]^2=k^3\cdot G(k,t). \]
Then from the condition $G(k,t)=0$ we obtain a non-trivial solution
of (\ref{1}):
\begin{equation} \label{2}
k_2=\frac{16t(t+1)(2t^2+4t+1)}{16t^4+64t^3+76t^2+24t-1} \,.
\end{equation}
Using (\ref{2}) we obtain the following expressions for $b$, $d$ and $e$:
\begin{equation} \label{3}
b(t)=\frac{16t(t+1)(t+2)(2t^2+4t+1)}{16t^4+64t^3+76t^2+24t-1} \,,
\end{equation}
{\footnotesize
\begin{equation} \label{4}
d(t)=\frac{256t^8+2048t^7+6272t^6+8960t^5+5424t^4+192t^3
 -888t^2-112t+33}
{16(16t^4+64t^3+76t^2+24t-1)(2t^2+4t+1)(t+1)}  \,,
\end{equation} }

{\footnotesize
\begin{equation} \label{5}
  \begin{array}{c}
 e(t)=
 (4096t^{12}+49152t^{11}+262144t^{10}+819200t^9+1665024t^8
+2310144t^7 \vspace{1ex}\\
 \mbox{}+2233728t^6+1507584t^5
+697856t^4+211968t^3+38624t^2+3520t+105) \vspace{1ex}\\
\mbox{} /[16(16t^4+64t^3+76t^2+24t-1)(2t^2+4t+1)(t+1)]  \,.
  \end{array}
\end{equation} }

\begin{theorem} \label{t:1}
Let $b(t)$, $d(t)$ and $e(t)$ be defined by (\ref{3}), (\ref{4}) and
(\ref{5}). Then the elliptic curve
\begin{equation} \label{6}
E: \quad y^2=(b(t)x+1) (d(t)x+1) (e(t)x+1)
\end{equation}
over $\Q(t)$ has the torsion group isomorphic to $\Z/2\Z \times \Z/2\Z$ and
the rank greater or equal $4$.
\end{theorem}

\dz
The points ${\cal O}$, $A=(-\frac{1}{b(t)},0)$, $B=(-\frac{1}{d(t)},0)$,
$C=(-\frac{1}{e(t)},0)$ form a subgroup of the torsion group
$E_{\rm tors}(\Q(t))$ which is isomorphic to $\Z/2\Z \times \Z/2\Z$.
By Mazur's theorem and a theorem of Silverman 
(\cite[Theorem 11.4, p.271]{Sil}), 
it suffices to check that there is no point on
$E(\Q(t))$ of order four or six.

If there is a point $D$ on $E(\Q(t))$ such that $2D\in \{A,B,C\}$, then
2-descent Proposition (see \cite[4.1, p.37]{Hus}) implies that
at least one of the expressions
$\pm b(t)[e(t)-d(t)]$, $\pm d(t)[e(t)-b(t)]$, $\pm e(t)[d(t)-b(t)]$ is a
perfect square. However, by specialization $t=1$ we see that this is not
the case.

If there is a point $F=(x,y)$ on $E(\Q(t))$ such that $3F=A$, $F\neq A$, then
from $2F=-F+A$ we obtain the equation
\begin{equation} \label{z}
 x^4 -6g(t)x^2-4g(t)h(t)x-3h(t)^2=0 \,,
\end{equation}
where $g(t)=b(t)e(t)+d(t)e(t)-2b(t)d(t)$,
$h(t)=b(t)d(t)[e(t)-d(t)][e(t)-b(t)]$. One can easily check that e.g.
for $t=1$ the equation
(\ref{z}) has no rational solution. Similarly we can prove that
there is no point $F$ on $E(\Q(t))$ such that $3F=B$, $F\neq B$ or
$3F=C$, $F\neq C$. Therefore, we conclude that
$E_{\rm tors}(\Q(t))$ is isomorphic to $\Z/2\Z \times \Z/2\Z$.

Now, we will prove that four points with $x$-coordinates $0$,
\[ a(t)=\frac{16t^4+64t^3+76t^2+24t-1}{16(2t^2+4t+1)(t+1)} \,,\]
\[ c(t)=\frac{(t+1)(16t^4\!+64t^3\!+68t^2\!+8t+1)
   (16t^4\!+64t^3\!+100t^2\!+72t+17)}
   {4(16t^4+64t^3+76t^2+24t-1)(2t^2+8t+1)} \]
and
\[ \frac{1}{b(t)d(t)e(t)} \]
are independent $\Q(t)$-rational points. Since the specification map
is always a homomorphism, we only have to show that there is a
rational number $t$ for which above four points are specialized to four
independent $\Q$-rational points. We claim that this is the case
for $t=1$.

We obtain the elliptic curve
\[ \begin{array}{c}
E^{*}: \quad
y^2 = x^3+6039621860663185x^2 \\
  \mbox{}+4139229575576935297875399628800x \\
  \mbox{}+48358738060886226093564403421659325399040000
\end{array}
\]
and the points
\[ P=(0,6954044726695840435200), \]
\[ Q=(2322788497348275,234053443113019268212650), \]
\[ R=(48986399479921200,11499867835919119918338000), \]
\[ S=(51511970169856/9,229496624258539337814016/27). \]
Then $S=2S_1$, $P-Q=2Q_1$, $P-R=2R_1$, where
\[ Q_1=(265264199014080, -39874704566573066299200), \]
\[ R_1=(3714953903426304, 387359212888080790925568), \]
\[ S_1=(2452641432447360, 247558457515476853468800). \]
It is sufficient to prove that the points $P,Q_1,R_1,S_1$ are
independent. The curve $E^{*}$ has three $2$-torsion points:
\[ A^{*}=(-11888861752320,0), \]
\[ B^{*}=(-5253470166461440,0), \]
\[ C^{*}=(-774262832449425,0). \]
Consider all points of the form
\[ X={\varepsilon}_1P+{\varepsilon}_2Q_1+{\varepsilon}_3R_1+
   {\varepsilon}_4S_1+T, \]
where ${\varepsilon}_i \in \{0,1\}$ for $i=1,2,3,4$, $T\in \{{\cal O},
A^{*},B^{*},C^{*}\}$ and $X=(x,y)\neq {\cal O}$. For all of these $63$ points
at least one of the numbers
$x+11888861752320$ and $x+5253470166461440$ is not a perfect square.
Hence, from 2-descent Proposition (\cite[4.1, p.37]{Hus}) it
follows that $X\not\in 2E(\Q)$.

Assume now that $P,Q_1,R_1,S_1$ are dependent modulo torsion, i.e.
that there exist integers $i,j,m,n$ such that
$|i|+|j|+|m|+|n|\neq 0$ and
\[ iP+jQ_1+mR_1+nS_1=T, \]
where $T\in \{{\cal O},A^{*},B^{*},C^{*}\}$.
Then the result which we just proved
shows that $i,j,m,n$ are even, say $i=2i_1$, $j=2j_1$, $m=2m_1$,
$n=2n_1$, and $T={\cal O}$. Thus we obtain
\[ i_1P+j_1Q_1+m_1R_1+n_1S_1 \in \{{\cal O},A^{*},B^{*},C^{*}\}. \]
Arguing as above, we conclude that $i_1,j_1,m_1,n_1$ are even,
and continuing this process we finally obtain that
$i=j=m=n=0$, a contradiction.
\qed

\medskip

By a theorem of Silverman (\cite[Theorem 11.4, p. 271]{Sil}),
the specialization map is an injective homomorphism for all but
finitely many points $t\in \Q$. This fact implies that by specialization
of the parameter $t$ to a rational number one gets in all but
finitely many cases elliptic curves over $\Q$ of rank at least four, 
and with subgroup of the torsion group which is 
isomorphic to $\Z/2\Z \times \Z/2\Z$. Hence, we have

\begin{corollary} \label{k:2}
There is an infinite number of elliptic curves over $\Q$ with 
three non-trivial 2-torsion points whose rank is greater or equal $4$.
\end{corollary}

\section{An example of high-rank curve}
We used the program {\it mwrank} (see \cite{Cre}) for computing the rank
of elliptic curves obtained from (\ref{6}) by specialization of
parameter $t$. However, since the coefficients in the corresponding
Weierstra\ss{} form are usually very large, we were able to determine the rank
unconditionally only for a few values of $t$. The following table shows the
values of $t$ for which we were able to compute the rank.

\bigskip

\begin{center}
\begin{tabular}{|c|cccccccccccc|} \hline
&&&&&&&&&&&&\\
$t$ &  \mbox{} &\mbox{} & $\frac{1}{4}$ & \mbox{} & $\frac{1}{2}$ & \mbox{}
& $1$ & \mbox{} &$\frac{3}{2}$& \mbox{} &$2$ &\\
&&&&&&&&&&&& \\ \hline
&&&&&&&&&&&&\\
Selmer rank & \mbox{} & \mbox{} & $8$ & \mbox{} & $4$ & \mbox{} & $5$ &
\mbox{} & $8$ & \mbox{} & $9$& \\
&&&&&&&&&&&& \\\hline
&&&&&&&&&&&&\\
rank &\mbox{} & \mbox{}  &$4$& \mbox{} & $4$& \mbox{}
&$5$ & \mbox{} & $6$&\mbox{} & $7$& \\
&&&&&&&&&&&& \\ \hline
\end{tabular}
\end{center}

\bigskip

Hence, we obtain

\begin{theorem} \label{t:3}
\,\,There is an elliptic curve over \,$\Q$\, with the torsion group
\linebreak $\Z/2\Z \times \Z/2\Z$ whose rank is equal 7.
\end{theorem}

Let us write this example of the curve with rank equal $7$ explicitly:
\begin{eqnarray*}
y^2=(\frac{2176}{373}x+1)(\frac{192386145}{101456}x+1)
(\frac{122265}{101456}x+1)
\end{eqnarray*}
or in Weierstra\ss{} form:
\begin{equation} \label{7}
  \begin{array}{c}
y^2= x^3+19125010376436745905x^2 \\
\mbox{}+52038165131253677052054066913723699200x \\
\mbox{}+521987941186440643611574160434960523120404754595840000 \,.
  \end{array}
\end{equation}
Seven independent points on (\ref{7}) are
\[ P_1=(727040606274688800,6989234854370183719797420000), \]
{\footnotesize
\[ P_2=(\frac{106210585076366036700000}{12769},
      \frac{69679298576214445317616490513378400}{1442897} ), \] }
{\footnotesize
\[ P_3=(\frac{335675366319765814629760}{71289},
      \frac{529539341511970538352844395949129600}{19034163} ), \] }
{\footnotesize
\[ P_4=(\frac{8891873190221412964144}{81},
      \frac{910251624041798036784012061900208}{729} ) ,\] }
{\footnotesize
\[ P_5=(\frac{101700294221755145291440}{841},
      \frac{34956857441184030025736520646806800}{24389} ) ,\] }
{\footnotesize
\[ P_6=(\frac{73133606420424854742955}{114921},
      \frac{251397104609526457099162042379450150}{38958219} ) ,\] }
\[ P_7=(-11146430015095060400,20291973801839968429609236400)  \vspace{3ex}.\]

{\large\bf Acknowledgements}.  The author would like to thank
J.E. Cremona and the referee for helpful suggestions and L. Kulesz
for kindly providing his recent preprint.

{\footnotesize
\nopagebreak{\sc Department of Mathematics,
University of Zagreb, Bijeni\v cka cesta 30,
10000 Zagreb, Croatia}

{\em E-mail address}: {\tt duje@math.hr} }

\end{document}